\documentclass[12pt]{article}

\usepackage[T2A]{fontenc}
\usepackage[cp866]{inputenc}

\usepackage{amsfonts}
\usepackage{amssymb}
\usepackage{amsmath}
\usepackage{amsthm}

\def \le {\leqslant}
\def \ge {\geqslant}

\begin{document}

\begin{Large}
\centerline{\bf Badly approximable numbers related to the Littlewood conjecture}
\end{Large}
 \vskip+1.5cm \centerline{\bf Moshchevitin N.G. \footnote{ Research is supported by
grants RFFI 06-01-00518  }} \vskip+1.5cm

\begin{small}
\centerline{\bf Abstract.}
By means of Peres-Schlag's method we prove the existence of real numbers $\alpha, \beta$ such that
$$
\liminf_{ q\to \infty} (q\log^2 q)||\alpha q|| \,||\beta q|| > 0.
$$
\end{small}

 \vskip+1.5cm

{\bf 1.} The famous Littlewood conjecture suppose that for any two real numbers  $\alpha, \beta $ one has
$$
\liminf_{ q\to \infty} q||\alpha q|| \,||\beta q|| =0
$$
(here $||\cdot ||$ denotes the distance to the nearest integer).

This conjecture is obviously true for {\it almost all} (in the
sense of the Lebesgue measure) pairs $(\alpha ,\beta ) \in
\mathbb{R}^2$. Moreover, from Gallagher's theorem (see \cite{KG})
we know that for almost all $\alpha, \beta $ one has
$$
\liminf_{ q\to \infty} (q\log^2 q)||\alpha q|| \,||\beta q|| =0.
$$
Here we prove the following result.

{\bf Theorem 1.}\,\,\,{\it There exist real numbers $\alpha, \beta$ such that
$$
\liminf_{ q\to \infty} (q\log^2 q)||\alpha q|| \,||\beta q|| > 0.
$$}
Our result is based on an original method of constructing badly approximable numbers invented by Y. Peres and W. Schlag in \cite{P}. This method
was recently applied by the author to several Diophantine problems \cite{M1}-\cite{M4}. In the sequel we suppose that we deal with the dyadic
logarithm $\log x = \log_2 x$.

{\bf 2.   }\,\,\, In this section we prove an elementary lemma dealing with the  sum
$$
\sigma (\alpha; q,Q) = \sum_{x=q}^{Q-1} \frac{1}{||\alpha x||\, x\log_2^2 x}.
$$
   We take $\alpha \in \mathbb{R}$ to be a badly
approximable number. This means that with some positive $\delta <1$ one has
\begin{equation}
\inf_{ q\in \mathbb{N}} q||\alpha q|| \ge \delta. \label{1}
\end{equation}

{\bf  Lemma 1. } \,\, {\it Let (\ref{1}) holds. Then for any natural numbers $Q>q$ one has
$$
\sigma (\alpha; q,Q)\le 2^8\delta^{-1}\left(\log_2\left(\frac{\log_2Q}{\log_2 q}\right)+\log_2(1/\delta )\right).
$$
}

Proof. For natural $\mu, \nu $ define
$$
A_{\mu,\nu} = \{ x\in \mathbb{N}:\,\,\, 2^\nu \le x < 2^{\nu+1},\,\, 2^{-\mu -1}< ||\alpha x||\le 2^{-\mu }\} .
$$
As for $2^\nu \le x < 2^{\nu+1}$  from (\ref{1}) follows
$$
\delta 2^{-\nu-1} <||\alpha x||\le 2^{-1},
$$
we see that
$$
\{x\in \mathbb{N}:\,\,\, 2^\nu \le x < 2^{\nu+1}\} = \bigsqcup_{\mu =1}^{\nu+1+\lceil\log_2 (1/\delta )\rceil } \,\,\, A_{\mu,\nu} . $$
   Now
$$
\sigma (\alpha; q,Q) \le \sum_{\nu= \lfloor\log_2 q\rfloor }^{\lceil\log_2 Q \rceil} \sum_{x=2^\nu}^{2^{\nu+1}-1} \frac{1}{||\alpha x||\,
x\log_2^2 x}=
  \sum_{\nu= \lfloor\log_2 q\rfloor }^{\lceil\log_2 Q \rceil} \sum_{\mu= 1 }^{\nu+1+\lceil\log_2 (1/\delta )\rceil}
\sum_{x\in A_{\mu,\nu} }\frac{1}{||\alpha x||\, x\log_2^2 x} \le$$
$$
\le \sum_{\nu= \lfloor\log_2 q\rfloor }^{\lceil\log_2 Q \rceil} \sum_{\mu= 1 }^{\nu+1+\lceil\log_2 (1/\delta )\rceil} 2^{\mu-\nu+1}\nu^{-2} {\rm
card}(A_{\mu,\nu}).
$$
Note that
$$
{\rm card}(A_{\mu,\nu})\le {\rm card}(\{x\in \mathbb{N}:\,\,\, x < 2^{\nu+1},\,\,||\alpha x||\le 2^{-\mu }\})=I_{\mu,\nu}.
$$
We have
\begin{equation}
\sigma (\alpha; q,Q)\le   \sum_{\nu= \lfloor\log_2 q\rfloor }^{\lceil\log_2 Q \rceil} \sum_{\mu= 1 }^{\nu+1+\lceil\log_2 (1/\delta )\rceil}
2^{\mu-\nu+1}\nu^{-2}  I_{\mu,\nu} . \label{I}
\end{equation}

 The value of $I_{\mu,\nu}$  is bounded by the number of integer points in the region
$$
\Omega  =
 \left\{ (x,y ) \in \mathbb{R}^2:\,\,\, 0\le x < 2^{\nu+1},\,\, |x \alpha  -y|\le 2^{-\mu} \right\}.
$$
Obviously the Lebesgue  measure of $\Omega $ is equal to $ {\rm mes} \Omega = 2^{\nu-\mu +1}$.

Let there exists an integer primitive point  $(p_0,a_0)\in \Omega
 $ (otherwise $ I_{\mu,\nu}=0$).  From (\ref{1}) we deduce
$$
\delta p_0^{-1}\le ||p_0\alpha ||\le 2^{-\mu}
$$
and $ p_0 \ge \delta 2^\mu$. Now we consider two cases.

In the {\it first case} we suppose that  $(p_0,a_0)$ is the unique primitive  integer point in $\Omega (p)$.
 Then
 $$
 I_{\mu,\nu }\le
 \lfloor2^{\nu+1}/p_0\rfloor+1\le
2^{\nu-\mu+2}\delta^{-1}.
$$

In the {\it second case} the convex hull of all integer points in the region $\Omega  )$ is a convex polygon $\Pi$ with integer vertices.
According to Pick's formula the number $I_{\mu,\nu }$ of integer point in the polygon $\Pi$ is not greater than
 $6{\rm mes}\Pi
$ and so
$$
I_{\mu,\nu }\le 6{\rm mes}\Pi
  \le 6{\rm mes} \Omega =
  12 \cdot 2^{\nu -\mu}.
  $$
  In any case
  \begin{equation}
 I_{\mu,\nu }\le 2^{\nu-\mu+4}\delta^{-1}.
 \label{II}\end{equation}
Now we substitute (\ref{II}) into (\ref{I}) and obtain
$$
\sigma (\alpha; q,Q)\le 2^5 \delta{-1}\cdot \sum_{\nu= \lfloor\log_2 q\rfloor }^{\lceil\log_2 Q \rceil}
 (\nu+1+\lceil\log_2 (1/\delta )\rceil )  \nu^{-2}\le
2^6\delta{-1}\left( \sum_{\nu= \lfloor\log_2 q\rfloor }^{\lceil\log_2 Q \rceil}\nu^{-1} +\log_2(1/\delta )\right)\le
$$
$$
\le 2^8\delta{-1}\left(\log_2\left(\frac{\log_2Q}{\log_2 q}\right)+\log_2(1/\delta )\right)
$$
and Lemma 1 is proved.

{\bf  3.  } In this section we (following Y. Peres and W. Schlag \cite{P}) construct "dangerous" sets of reals. Let $\varepsilon$ be positive
and small enough.
 For integers $ 2\le x, 0\le y\le
x$  define
\begin{equation}
E (x,y) = \left[ \frac{y}{x} -\frac{\varepsilon}{||\alpha x || x^{2}\log_2^2 x}, \frac{y}{x} +\frac{\varepsilon}{||\alpha x || x^{2}\log_2^2 x}
\right],
 \,\,\,
E (x) =\bigcup_{y =0}^{x} E (x,y)\bigcap [0,1]. \label{E}
\end{equation}
 Define
\begin{equation}
l_0 = 0,\,\,\,
 l_x = \lfloor\log_2 (||\alpha x || x^{2}\log_2^2 x /2\varepsilon )  \rfloor ,\,\, x \in \mathbb{N} . \label{L}
\end{equation}
 Each
segment form the union $E_\alpha (x)$ from (\ref{E}) can be covered by a dyadic interval of the form
$$
\left( \frac{b}{2^{l_x }}, \frac{b+z}{2^{l_x }}\right),\,\,\, z = 1,2 .$$

Let $A (x)$ be the smallest union of all such dyadic segments which cover the whole set  $E (x)$. Put
$$
A^c (x) = [0,1] \setminus A (x). $$ Then
$$
A^c (x)
 = \bigcup_{\nu = 1}^{\tau_x } I_\nu
$$
where closed  segments $I_\nu $ are of the form
\begin{equation}
\left[ \frac{a}{2^{l_x}}, \frac{a+1}{2^{l_x}}\right] ,\,\,\, a\in \mathbb{Z}. \label{aaa}
\end{equation}

 We take $q_0$ to be a large  positive integer. In order to prove Theorem 1 it is sufficient to show that for all $
q \ge q_0 $ the sets
$$
B_q =\bigcap_{x=q_0}^q A^c_{} (x)$$ are not empty. Indeed as the sets $B_q$ are closed and nested we see that there exists real $\beta$ such
that
$$
\beta \in \bigcap_{q\ge q_0} B_q .$$ One can see that the pair $\alpha, \beta $ satisfies the conclusion of Theorem 1.

In the next section we prove the following statement.

{\bf Lemma 2.}\,\,{\it Suppose that $\varepsilon $ is  small enough. Then for $q_0$ large enough and for any
$$q_1\ge q_0
,\,\,\, q_2 = q_1^3,\,\,\, q_3 =q_2^3
$$ the following holds. If
\begin{equation}
{\rm mes} B_{q_2} \ge {\rm mes} B_{q_1}/2>0 \label{m1}
\end{equation}
then
\begin{equation}
{\rm mes} B_{q_3} \ge {\rm mes} B_{q_2}/2>0. \label{m2}
\end{equation}
} Theorem 1 follows from Lemma 2 by induction as the base of the induction  obviously follows from the arguments of Lemma's proof.

{\bf  4.} \,\,\, In this section we prove Lemma 2. Here we also follow the arguments from the paper \cite{P} by Y. Peres and W. Schlag. First of
all we show that for $x\ge q^3$ where $q\ge q_0$ one has
\begin{equation}
 {\rm mes }\left( B_q  \bigcap A  (x) \right) \le \frac{2^4\varepsilon}{||\alpha x||\, x\log_2^2 x} \times {\rm mes} B_q.
\label{PP}
\end{equation}
Indeed as from (\ref{L}) it follows  that  $$l_x \le l= \lfloor \log_2( q^{2}\log_2^2 q /\varepsilon ) \rfloor ,\, \forall x\le q$$ we see that
 $B_q$ is a union
$$
B_q = \bigcup_{\nu = 1}^{T_q } J_\nu
$$
with $J_\nu$ of the form $$\left[ \frac{a}{2^{l}}, \frac{a+1}{2^{l}}\right] ,\,\,\, a\in \mathbb{Z}.
$$
Note that $A(x)$ consists of the segments of the form (\ref{aaa}) and for $ x\ge q^3 > 2^{l+1}
 $ (for $q_0$ large enough)   we see that each $J_\nu$ has at least two rational fractions of the form $\frac{y}{x}, \frac{y+1}{x}$ inside. So
 \begin{equation}
{\rm mes} (J_\nu \cap A(x)) \le     \frac{2^4\varepsilon}{||\alpha x||\, x\log_2^2 x} \times {\rm mes} J_\nu. \label{aaaa}\end{equation} Nof
(\ref{PP}) follows from (\ref{aaaa}) by summation over $ 1\le \nu\le T_q$.

To continue we observe that
$$
B_{q_3} = B_{q_2} \setminus
 \left(\bigcup_{x=q_2+1}^{q_3}  A (x) \right)   ,
$$
and hence
$$
{\rm mes}  B_{q_3} \ge {\rm mes}   B_{q_2}  - \sum_{x=q_2+1}^{q_3} {\rm mes} ( B_{q_2}\cap A (x) ).
$$
As
 $$
B_{q_2}\cap A (x)\subseteq B_{q_1}\cap A_{\alpha_i} (x) $$ we can apply (\ref{PP}) for every $ x$ from the interval $q_1^3\le q_2 < x\le q_3$:
$$
{\rm mes}(B_{q_2}\cap A (x)) \le {\rm mes}(B_{q_1}\cap A (x)) \le \frac{2^4\varepsilon}{||\alpha x||\, x\log_2^2 x} \times {\rm mes} B_{q_1} \le
\frac{2^5\varepsilon}{||\alpha x||\, x\log_2^2 x} \times {\rm mes} B_{q_2}
$$
(in the last inequality  we use the condition  (\ref{m1}) of Lemma 2). Now as $\frac{\log_2q_3}{\log_2 q_2}= 3$ the conclusion (\ref{m2}) of
Lemma 2 follows from Lemma 1 for $\varepsilon $ small enough:
$$
{\rm mes}  B_{q_3}\ge {\rm mes}   B_{q_2} \left( 1- 2^5\varepsilon \sigma (\alpha; q_2+1,q_3+1) \right)\ge
$$
$$
\ge
 {\rm mes}   B_{q_2} \left( 1- 2^{14}\varepsilon \left(\log_2\left(\frac{\log_2q_3}{\log_2 q_2}\right)+\log_2(1/\delta )\right) \right)\ge
{\rm mes}   B_{q_2}/2.
$$

\vskip+2.0cm

author: Nikolay Moshchevitin

\vskip+0.5cm

e-mail: moshchevitin@mech.math.msu.su, moshchevitin@rambler.ru

\end{document}